\documentclass[12pt]{article}
\usepackage[utf8]{inputenc}
\usepackage[T1]{fontenc}
\usepackage{amsmath, amssymb}
\usepackage{geometry}
\usepackage[hidelinks]{hyperref}
\title{ § 12.\\
    {\small\textit{English translation of Frobenius' and Stickelberger's}} \\[0.5em] 
    \MakeUppercase{On the Theory of Elliptic Functions}
    \\[0.5em] {\small \textit{Journal für die reine und angewandte Mathematik} 83, 175--179 (1877)}
    }
\author{Ferdinand Georg Frobenius with Ludwig Stickelberger \\ Zürich }
\date{March 10, 1877}

\begin{document}
\maketitle

\medskip

A remarkable formula communicated by Mr. \textit{Hermite} in a recently published note on elliptic functions (this Journal Vol.~82, p.~343) prompts us to point out the connections between several related formulas. 
To pass from the general equation with which we begin to the more specialised formula of Mr. \textit{Hermite}, and from there to the further specialised one that Mr. \textit{Kiepert}\footnote{\textit{Kiepert}, Actual execution of the integer multiplication of elliptic functions, this Journal Vol.~76, p.~21. There one finds the definition of the functions $\sigma(u)$ and $\wp(u)$ which Mr. \textit{Weierstrass} has introduced into the theory of elliptic functions, as well as a brief compilation of their most important properties. Regarding the formulas and theorems from this theory which we will use in the following, we refer to that treatise.} has used as the basis of his solution to the multiplication problem, we employ a limiting process which we first present in its general form.

Let
\begin{equation*}
f_0(u),\quad f_1(u),\quad \dots \quad f_n(u)
\end{equation*}
be $n+1$ convergent series of positive powers of $u$, and let
$u_0,\quad u_1,\quad \dots \quad u_n$ be arbitrary values within their common region of convergence. Then the determinant of $(n+1)\textsuperscript{th}$ degree
\begin{equation*}
|f_\alpha(u_\beta)| = F(u_0, u_1, \dots u_n)
\end{equation*}
can be expanded into a series of powers of $u_0, u_1, \dots u_n$ and this, as an alternating function, can be brought to the form
\begin{equation*}
F = G(u_0, u_1, \dots u_n) \Pi(u_\alpha - u_\beta)
\end{equation*}
where the symmetric function $G$ is likewise a series of \textit{positive} powers of $u_0, u_1, \dots u_n$. 
(In the difference product, $\alpha$ and $\beta$ are to run through the pairs of numbers $0, 1, \dots n$ here and in the following such that $\alpha > \beta$.) Then the value which the function $G(u_0, u_1, \dots u_n)$ assumes when all
$n+1$ variables are set equal to $u$ can easily be expressed as a determinant. For simplicity, let $h$ denote a small quantity and set
\begin{equation*}
u_\beta = u + \beta h, \quad (\beta = 0, 1, \dots n)
\end{equation*}
and introduce the notation
\begin{equation*}
\Delta f(u) = f(u+h) - f(u),
\end{equation*}
then by a known determinant theorem
\begin{equation*}
|f_\alpha(u_\beta)| = |\Delta^\beta f_\alpha(u)|
\end{equation*}
and consequently
\begin{equation*}
G = \frac{F}{\Pi(u_\alpha - u_\beta)} = \frac{1}{\Pi(\alpha-\beta)} \left| \frac{\Delta^\beta f_\alpha(u)}{h^\beta} \right|.
\end{equation*}

From this one obtains, by letting $h$ approach the limit $0$, the sought relation
\begin{equation}
\lim \frac{F(u_0, u_1, \dots u_n)}{\Pi(u_\alpha-u_\beta)} = \frac{1}{\Pi(\alpha-\beta)} |f_\alpha^{(\beta)}(u)|.
\end{equation}

On the left side, the quotient is first to be expanded into a series of positive powers of $u_0, u_1, \dots u_n$, and then the arguments are all to be set equal to $u$.

We now apply a similar procedure to the determinant
\begin{equation*}
R = \begin{vmatrix} 0 & 1 & \dots & 1 \\ 1 & \psi(u_0+v_0) & \dots & \psi(u_0+v_n) \\ . & . & \dots & . \\ 1 & \psi(u_n+v_0) & \dots & \psi(u_n+v_n) \end{vmatrix},
\end{equation*}
where
\begin{equation*}
\psi(u) = \frac{d\log\sigma(u)}{du}.
\end{equation*}
The difference $\psi(u+v) - \psi(u)$ is a doubly periodic function of $u$. By multiplying the elements of the first row of $R$ with $\psi(u_0)$ and subtracting them from those of the second row, one recognises that this determinant is a doubly periodic function of $u_0$. Since the same conclusion applies to the other variables entering into $R$, this function of $2n+2$ arguments also has the remarkable property of being doubly periodic with respect to each of them.

The function $\psi(u)$ becomes infinite only for $u = 0$ (and congruent values), and its expansion in ascending powers of $u$ begins with $\frac{1}{u}$.
Considered as a function of $u_0$, $R$ therefore becomes infinite only at the $n+1$ values
\begin{equation*}
u_0 = -v_0, \quad -v_1, \quad \dots \quad -v_n
\end{equation*}
and congruent values, with simple poles at each. On the other hand, $R$ obviously vanishes for
\begin{equation*}
u_0 = u_1, \quad \dots \quad u_n.
\end{equation*}

Now however an elliptic function\footnote{Following Mr. \textit{Weierstrass}, we call a doubly periodic function elliptic if it has the character of a rational function everywhere in the finite plane, i.e.\ if it can be expanded in the neighborhood of every finite value $a$ as a series in integer powers of $u-a$ containing only finitely many negative powers.} has as many zeros as infinities, and (after Abel's theorem) the sum of the values for which it vanishes is congruent to the sum of the values for which it becomes infinite. (\textit{Briot et Bouquet}, Fonctions elliptiques, II.~éd., p.~241, Théor.~III; p.~242, Théor.~V. \textit{Kiepert}, l.~c.~p.~24 and 25.) Consequently $R$ must also vanish for an $(n+1)$th value of $u_0$, which is to be calculated from the equation
\begin{equation*}
u_0 + v_0 + \dots + u_n + v_n = 0
\end{equation*}
By a known theorem (\textit{Briot et Bouquet}, p.~242, Théor.~IV; p.~243, Théor.~VI. \textit{Kiepert}, l.~c.) that determinant is therefore, up to a factor independent of $u_0$, equal to
\begin{equation*}
\frac{\sigma(u_0+v_0+\dots+u_n+v_n)\sigma(u_1-u_0)\sigma(u_2-u_0)\dots\sigma(u_n-u_0)}{\sigma(u_0+v_0)\sigma(u_0+v_1)\dots\sigma(u_0+v_n)}.
\end{equation*}
If one investigates in a similar way its dependence on the remaining $2n+1$ arguments, one finds\footnote{For $n=1$ this formula essentially coincides with that which \textit{Jacobi} has given in Vol.~15 of this Journal (p.~204, 13).}, up to a constant factor, that,
\begin{equation}
\begin{aligned}
-\begin{vmatrix}
0 & 1 & \dots & 1 \\
1 & \frac{\sigma'(u_0+v_0)}{\sigma(u_0+v_0)} & \dots & \frac{\sigma'(u_0+v_n)}{\sigma(u_0+v_n)} \\
. & . & \dots & . \\
1 & \frac{\sigma'(u_n+v_0)}{\sigma(u_n+v_0)} & \dots & \frac{\sigma'(u_n+v_n)}{\sigma(u_n+v_n)}
\end{vmatrix} 
= \frac{\sigma(u_0+v_0+\dots+u_n+v_n)\Pi\sigma(u_\alpha-u_\beta)\Pi\sigma(v_\alpha-v_\beta)}{\Pi\sigma(u_\alpha+v_\beta)}.
\end{aligned}
\end{equation}
In the denominator of the right side the product is to be extended over all pairs of numbers from $0$ to $n$, and in the numerator only over those for which $\alpha > \beta$.

To verify that the constant factor in equation (2.) is correctly given, one observes that it holds immediately for $n=0$, and in general follows easily by induction from $n$ to $n+1$. Indeed, if one multiplies the elements of the last row on the left side by $u_n+v_n$ and then sets $u_n = -v_n$, they all vanish except the last, which becomes $1$. Therefore $R$ reduces to the determinant formed analogously from the $n$ argument pairs $u_0, v_0, \dots u_{n-1}, v_{n-1}$. In the expression on the right side of equation (2.), when $u_n = -v_n$ we have
\begin{equation*}
\lim \frac{u_n+v_n}{\sigma(u_n+v_n)} = 1,
\end{equation*}
\begin{equation*}
\frac{\sigma(u_n-u_\beta)\sigma(v_n-v_\beta)}{\sigma(v_n+u_\beta)\sigma(u_n+v_\beta)} = 1 \quad (\beta = 0, 1, \dots n-1).
\end{equation*}
Thus the right side likewise reduces to the expression formed analogously from the $n$ argument pairs $u_0, v_0, \dots u_{n-1}, v_{n-1}$.

In equation (2.) we now set
\begin{equation*}
v_\beta = \beta h \quad (\beta = 0, 1, \dots n)
\end{equation*}
and transform the left side by the method already applied above into a determinant in which the $(\alpha+2)\textsuperscript{th}$ row contains the elements
\begin{equation*}
1, \quad \psi(u_\alpha), \quad \Delta\psi(u_\alpha), \quad \Delta^2\psi(u_\alpha), \quad \dots \quad \Delta^n\psi(u_\alpha).
\end{equation*}
Since in this transformation the elements of the first row become
\begin{equation*}
0, \quad 1, \quad 0, \quad 0, \quad \dots \quad 0
\end{equation*}
the left side reduces to a determinant of $(n+1)\textsuperscript{th}$ degree, which we denote in an easily understandable manner with
\begin{equation*}
-R = |1, \quad \Delta\psi(u_\alpha), \quad \Delta^2\psi(u_\alpha), \quad \dots \quad \Delta^n\psi(u_\alpha)|.
\end{equation*}
Consequently
\begin{equation*}
\left| 1, \quad \frac{\Delta\psi(u_\alpha)}{h}, \quad \frac{\Delta^2\psi(u_\alpha)}{h^2}, \quad \dots \quad \frac{\Delta^n\psi(u_\alpha)}{h^n} \right|
\end{equation*}
\begin{equation*}
= \frac{\sigma(u_0+v_0+\dots+u_n+v_n)\Pi\sigma(u_\alpha-u_\beta)}{\Pi\sigma(u_\alpha+v_\beta)} \Pi\frac{\sigma(v_\alpha-v_\beta)}{h},
\end{equation*}
therefore, in the limit that $h$ approaches zero,
\begin{equation*}
|1, \quad \psi'(u_\alpha), \quad \psi''(u_\alpha), \quad \dots \quad \psi^{(n)}(u_\alpha)| = \frac{\sigma(u_0+\dots+u_n)\Pi\sigma(u_\alpha-u_\beta)\Pi(\alpha-\beta)}{(\Pi\sigma(u_\alpha))^{n+1}}.
\end{equation*}
Setting according to Mr. \textit{Weierstrass}
\begin{equation*}
-\psi'(u) = -\frac{d^2\log\sigma(u)}{du^2} = \wp(u),
\end{equation*}
this formula reads
\begin{equation}
\begin{vmatrix}
1 & \wp(u_0) & \wp'(u_0) & \dots & \wp^{(n-1)}(u_0) \\
1 & \wp(u_1) & \wp'(u_1) & \dots & \wp^{(n-1)}(u_1) \\
. & . & . & \dots & . \\
1 & \wp(u_n) & \wp'(u_n) & \dots & \wp^{(n-1)}(u_n)
\end{vmatrix}
= \frac{(-1)^n\Pi(\alpha-\beta)\sigma(u_0+\dots+u_n)\Pi\sigma(u_\alpha-u_\beta)}{(\Pi\sigma(u_\alpha))^{n+1}}.
\end{equation}

\medskip

This is essentially the equation given by Mr. \textit{Hermite} in the above cited letter. We now apply the limiting process from the beginning once more, choosing in formula (1.)
\begin{equation*}
f_0(u) = 1, \quad f_1(u) = \wp(u), \quad \dots \quad f_n(u) = \wp^{(n-1)}(u).
\end{equation*}
Then in the determinant $|f_\alpha^{(\beta)}(u)|$ the elements of the first column vanish except for the first, and it reduces to a determinant of $n\textsuperscript{th}$ degree. One thus arrives at the formula of Mr. \textit{Kiepert} (l.~c.~p.~31, formula ($29^a$.))
\begin{equation}
\begin{vmatrix}
\wp'(u) & \wp''(u) & \dots & \wp^{(n)}(u) \\
\wp''(u) & \wp'''(u) & \dots & \wp^{(n+1)}(u) \\
. & . & \dots & . \\
\wp^{(n)}(u) & \wp^{(n+1)}(u) & \dots & \wp^{(2n-1)}(u)
\end{vmatrix}
= \frac{(-1)^n(\Pi(\alpha-\beta))^2\sigma((n+1)u)}{\sigma(u)^{(n+1)(n+1)}}.
\end{equation}

\vspace{1in}
\begin{center}
\raisebox{0.5ex}{\rule{1in}{0.5pt}}\quad$\bullet$\quad$\diamond$\quad$\bullet$\quad\raisebox{0.5ex}{\rule{1in}{0.5pt}}
\end{center}

\newpage
\section*{Translator's Notes}

\vspace{1em}
\noindent\textbf{Original paper.} Frobenius, F. G., \& Stickelberger, L. (1877). 
Zur Theorie der elliptischen Functionen. 
\textit{Journal für die reine und angewandte Mathematik}, 83, 175--179. \\[0.5em]
DOI: \url{https://doi.org/10.1515/crll.1877.83.175} \\
Wikimedia scans: 

{\footnotesize\url{https://en.wikipedia.org/wiki/File:Zur_Theorie_der_elliptischen_Functionen.djvu}}
\vspace{1em}

\noindent{Translated by K. Khazanzheiev and G. Hesketh, March 2026.}

\end{document}